\documentclass[12pt]{article}%%%
\usepackage{setspace}%%%
\usepackage[vmargin=2.5cm, hmargin=3.0cm]{geometry}%%%
\usepackage{color} % for colored text%%%
\definecolor{mygray}{gray}{0.25} % mygray (0 black through 1 white)%%%
\usepackage{amssymb} \usepackage{amsmath,amscd,mathrsfs}

\usepackage{graphpap}
\usepackage{pstricks}
\usepackage{pst-all}
\usepackage{pstcol}
\usepackage{color}
\usepackage{enumitem}
\usepackage{hyperref}
\usepackage[active]{srcltx}

\newtheorem{thm}{Theorem}
\newtheorem{defin}{Definition}
 
\newtheorem{lem}{Lemma}

\newtheorem{exam}{Example}

\def\C{{ \! \rm \ I\!\!\!C}}

\def\N{{ \! \rm \ I\!N}}

\def\R{{ \! \rm \ I\!R}}

\newcommand{\parcial}[2]{\frac{\partial#1}{\partial#2}}
\newcommand{\norm}[1]{\vert \vert #1 \vert \vert}
\newcommand{\normv}[1]{ \vert #1 \vert }

\title{Absence of singular continuous spectrum for some geometric Laplacians} \author{Leonardo A. Cano Garc\'{i}a}
\begin{document}
\maketitle
\onehalfspacing%%%
\begin{abstract}
We provide two examples of spectral analysis techniques of
Schr\"odinger operators applied  to  geometric Laplacians. In
particular we  show how to adapt the method of analytic dilation
 to Laplacians on complete  manifolds with corners of
codimension 2 finding the absence of singular continuous spectrum
for these operators, a description of the behavior of its pure
point spectrum in terms of the underlying geometry, and a theory
of quantum resonances.
\end{abstract}
\section*{Introduction}
Spectral geometry  studies the interactions between the geometry
of a Riemannian manifold and the spectral analysis of its associated
Laplacian. These interactions have been deeply studied in the case
in which the manifold is closed (see \cite{Chavel} \cite{Rosemberg}). In
the closed case the Laplacian is a self-adjoint operator with compact
resolvent and hence its spectrum is purely discrete. This contrasts
with the general case since, if $(M,g)$ is a geodesically complete
Riemannian manifold, not necessarily closed, it is known that its
Laplacian $\Delta_g:C^\infty_c(M) \to L^2(M,dvol_g)$ is essentially self-adjoint  but its spectral resolution is not purely discrete in general.

In order to show more clearly the new spectral phenomenons for
$\Delta_g$ that the lost of compactness of $M$ brings with, let us
remind the spectral theorem for self-adjoint operators.
\begin{thm}\label{thm:spectral thm}{\upshape \cite{RS1}} Let $A$ be
a self-adjoint operator acting on a separable Hilbert space
$\mathscr{H}$. Then, there exists
 a collection of measures  on $\R$, $\{\mu_i\}_{i \in I}$  where
  $I \subseteq \N$, and there exists a unitary operator $\mathscr{F}$ such that
  the following diagram commutes
 \[
   \begin{CD}
   \mathscr{H} @>{A}>> \mathscr{H}\\
   @A{\mathscr{F}}AA@VV{\mathscr{F}^*}V\\
   \bigoplus_{i \in I} L^2(\R,\mu_i)@>>{\oplus_{i \in I}}x_i> \bigoplus_{i \in I}
   L^2(\R,\mu_i),
    \end{CD}
 \]
where the operator $\oplus_{i \in I}x_i$ is the direct sum of the
multiplication operators $x_i:L^2(\R,\mu_i) \to L^2(\R,\mu_i)$
that, specifically, send a function $f \in L^2(\R,\mu_i)$ to the
function $x_i f$.
\end{thm}

Given a self-adjoint operator $A$, as in the  theorem above, the
Lebesgue decomposition theorem for measures induces a
decomposition of the Hilbert space $\mathscr{H}$ into three
important invariants subspaces $\mathscr{H}_{ac}$,
$\mathscr{H}_{pp}$ and $\mathscr{H}_{sing}$ that, using theorem
\ref{thm:spectral thm}, are the associated to the three  Hilbert
spaces $$\mathop{\bigoplus}_{\mu_i \text{ is ac}}L^2(\R,\mu_i),\,
\mathop{\bigoplus}_{\mu_i \text{ is pp}}L^2(\R,\mu_i) \text{ and
}\mathop{\bigoplus}_{\mu_i \text{ is s}}L^2(\R,\mu_i),$$ where ac
means absolutely continuous, pp means pure point and s means
singular continuous measures (all of them defined respect to the
Lebesgue measure in $\R$).  We have already said that for $(M,g)$
a closed Riemannian manifold, its Laplacian has purely discrete
spectrum, in other words the invariant subspaces
$L^2_{ac}(M,dvol_g)$ and $L^2_{sing}(M,dvol_g)$, corresponding to
the self-adjoint operator $A:=\Delta_g$,  are empty. In this paper
 we will show examples of manifolds $(M,g)$ for which
 $L^2_{ac}(M,dvol_g)$ is not  empty
 and $L^2_{sing}(M,dvol_g)$ is  empty.

The former description of the spectral analysis of self-adjoint
operators justifies the question about how to prove that a geometric
Laplacian does not have singular continuous spectrum.  This paper provide,
through the friendly environment of precise examples and without
pretending to give a complete answer, an illustration of how to
work around this question. The common feature of these
examples is the control on the Riemannian metric of the open
manifold at infinity.

Let us  describe more carefully the contents. For Laplacians
acting on manifolds with cylindrical or cusps ends and their
natural generalizations it is well known in the
literature (see \cite{GUILLOPE}\cite{HUS} \cite{MuellerCusp}) how
to find a meromorphic extension of their resolvent  and how, using
it,  to prove the absence of singular continuous spectrum. In section \ref{section:classic results}
we illustrate the classical method  to obtain a meromorphic
extension of the resolvent for manifolds with cylindrical and cusp
ends, and we show the relation of it with the absence of singular
 continuous spectrum, that is, fundamentally, the limit absorption principle.
In section \ref{sect:Mourreanddilat} we explain the notion of
complete manifold with corner of codimension 2 and how  to prove
absence of singular continuous spectrum for Laplacians on this kind of
manifolds using the method of analylic dilation. Section
\ref{sect:Mourreanddilat} describe following \cite{CANODILATION}  the way to apply the technique
of analytic dilation to Laplacians on complete manifolds with
corners of codimension 2, context on which the techniques
explained in section \ref{section:classic results} do not apply. The method of
analytic dilation was originally applied to $N$-particle
Schr\"odinger operators, a classic reference in that setting is
\cite{GERARD}. It has also been  applied to the black-box
perturbations of the Euclidean Laplacian in the series of papers
\cite{SjostrandZworski1} \cite{SjostrandZworski2}
\cite{SjostrandZworski3} \cite{SjostrandZworski4}. In
\cite{B} it is used to study Laplacians on hyperbolic manifolds.
The analytic dilation has also been applied to the study of the
spectral and scattering theory of quantum wave guides and
Dirichlet boundary domains, see e.g. \cite{DUEXMES} \cite{KOVSAC}.
It has also been applied to arbitrary symmetric spaces of
noncompact type in \cite{MV1} \cite{MV2} \cite{MV3}.
In \cite{KALVIN1} it is applied to manifolds with analytic
asimptotically cylindrical end. In each of these settings new
ideas and new methods carry out.
%------------------------------------------------------------------------------------------
%--------------------------------------------------------
\section{Meromorphic extension of the resolvent and singular continuous spectrum}
\label{section:classic results} In this section we give the main
ideas of a method for meromorphically extend the resolvent of
Laplacians on manifolds with   cylindrical and  cusps ends and
show why such extension is enough to have absence of singular
continuous spectrum. The results of this section were obtained in
\cite{GUILLOPE}\cite{MuellerCusp} but we base our exposition in
\cite{HUS} since we consider it is easier to understand for
non experts.
\begin{defin}\label{defin:manifold cyl end}
Let $M$ be an open manifold  with a decomposition in a compact
manifold $M_0$ with boundary $Y:=\partial M_0$ and an open
manifold $M_\infty$ with boundary and suppose that $\partial
M_\infty=M_0 \cap M_\infty=\partial M_0$. If $M$ is endowed with a
complete Riemannian metric such that
\begin{itemize}
\item[i)] $M_\infty$ is isometric to $Y \times \R_+$ with the natural
product metric $g_y+du \otimes du$, we say that {\bf $M$ is a
manifold with cylindrical end}. \item[ii)] $M_\infty$ is isometric to
$Y \times \R_+$ with the Riemannian metric $\frac{du \otimes
du}{u^2} +\frac{g_y}{u^2}$ , we say that {\bf $M$ is a a manifold
with cusp end}.
%\item $M_\infty$ is isometric to $\R^n-S$ where $S$ is some compact submanifold of $\R^n$ with boundary $Y$ then we say that {\bf $M$ is a a manifold with Euclidean end}.
\end{itemize}
\end{defin}
Let $\Delta_{cyl}$ be the Laplacian associated to a manifold $(M,g_{M})$ with
cylindrical end. We recall that on $M_\infty=Y \times
\R_+$ the Laplacian $\Delta_{cyl}$ has the form
\begin{equation*}
\Delta_{0,M}:=-\parcial{^2}{u^2}+\Delta_Y
 \end{equation*}
where $\Delta_Y$ is the Laplacian  on $Y$ associated to the
Riemannian metric $g_Y$.   Let $(N^n,g_{N})$ be a manifold with
cusp, with $N_\infty=Y \times \R_+$, and let $\Delta_{cusp}$ be
the Laplacian associated to $(N^n,g_{N})$. $\Delta_{cusp}$ has the
form
\begin{equation*}
\Delta_{0,N}:=-u^2 \parcial{^2}{u^2}+(n-1)u\parcial{}{u}+u^2
\Delta_Y
 \end{equation*}
on $N_\infty$. An important approach, in quantum scattering theory, is
to think that the Laplacians $\Delta_{cyl}$ and $\Delta_{cusp}$
are perturbations of the operators $\Delta_{0,M}$ and
$\Delta_{0,N}$ respectively. Roughly speaking, this idea is formalized in
the following way. Using Dirichlet boundary conditions, we
consider the operators $\Delta_{0,M}$ and $\Delta_{0,N}$ as
self-adjoint operators in $L^2(Y \times \R_+,dvol_{g_M})$ and
$L^2(Y \times \R_+,dvol_{g_N})$. Observe that $L^2(Y \times
\R_+,dvol_{g_M})$ is naturally isomorphic to $L^2(Y,dvol_{g_Y})
\otimes L^2(\R_+,du)$; furthermore, if $\phi_i \in C^\infty(Y)$ is
an orthonormal basis of eigenvectors of $\Delta_Y$ with
eigenvalues $\mu_i$, then it provides an isomorphism of
$L^2(Y,dvol_{g_Y}) \otimes L^2(\R_+,du)$ with $\oplus_{i \in \N}
L^2(\R_+,du)$, hence $L^2(Y \times \R_+,dvol_{g_M}) \cong
\oplus_{i \in \N} L^2(\R_+,du)$. Modulo this last isomorphism, we
have:
\begin{equation}\label{eq:Laplacian in cyl}
\Delta_{0,M}=\bigoplus_{i \in
\N}\left(-\parcial{^2}{u^2}+\mu_i\right).
 \end{equation}
The analogue of (\ref{eq:Laplacian in cyl})  in the context of
manifolds with cusps is
\begin{equation}\label{eq:Laplacian in cusp}
\Delta_{0,M}=\bigoplus_{i \in \N}\left(-u^2
\parcial{^2}{u^2}+(n-1)u\parcial{}{u}+u^2\mu_j\right).
 \end{equation}
From this point we continue our exposition over the manifold with
cylindrical end $M$, indicating how analogous methods apply in the
case of manifolds with cusps and refering to \cite{MuellerCusp}
for details.

Fourier transform and (\ref{eq:Laplacian in cusp}) give us a
formula for the resolvent of the operators
$-\parcial{^2}{u^2}+\mu_i$ and hence  for the resolvent of
$\Delta_{0,M}$. The analogue of such a formula in the case of the
manifold $N$ with cusp is technically harder to obtain and can be
found in  \cite[lemma 2.68]{MuellerCusp}. Define the double of
the compact manifold with boundary
$M_0$ as $\tilde{M}:=M_0 \cup_Y M_0$, where we are identifying
the boundary $Y$ of two disjoint copies of $M_0$ and we endow
$\tilde{M}$ with the natural differential and Riemmanian structure.
 We have also a nice formula
for the resolvent of  the Laplacian $\Delta_{\tilde{M}}$ of
$\tilde{M}$ using a spectral resolution of $\Delta_{\tilde{M}}$.
In order to apply this knowledge about the resolvents in
$M_\infty$ and $\tilde{M}$, we construct a {\it parametrix} for
the resolvent of $\Delta_M$ i. e. an operator $P(\lambda)$ such
that $R(\lambda)-P(\lambda)$ is compact in some weighted
$L^2$-space, where $\lambda \in \C-\R_+$ and
$R(\lambda):=(\Delta_M-\lambda)^{-1}$ denotes the resolvent of
$\Delta_M$. We proceed as in \cite{HUS}: for $0\geq a \geq b$, let
$\rho(a,b) \in C^\infty(M,[0,1])$ be such that
\begin{equation*}
\rho(a,b)(x)=\begin{cases}0 &\text{ for } x \in M \cup
(Y\times[0,a]);\\
1& \text{ for } x \in Y\times[b,\infty). \end{cases}
\end{equation*}
We define the functions:
\begin{equation*}
\begin{split}
&\Phi_1:=1-\rho(\frac{4}{5},1), \hspace{1.5cm} \Psi_1:=1-\rho(\frac{2}{5},\frac{3}{5})\\
&\Phi_2:=\rho(0,\frac{1}{5}), \hspace{2.3cm} \Psi_2:=1-\Psi_1,
\end{split}
\end{equation*}
for which $\Psi_1+\Psi_2=1$,
\begin{equation*}\label{eq:partition resolvent}
\begin{split}
 \Phi_j(x)=1 \text{ for }x \in supp(\Psi_j), \text{ and }
\text{dist}(supp \nabla \Phi_j,supp \Psi_j)\geq \frac{1}{5}.
\end{split}
\end{equation*}
We define the operator $S(\lambda)$ with Schwartz kernel
\begin{equation}\label{eq:kernel resolv}
S(x_1,x_2,\lambda):=\sum_{j=1}^2 \Psi_j(x_1) R_j(x_1,x_2, \lambda)
\Phi_j(x_2),
\end{equation}
where $R_1(x_1,x_2, \lambda)$ is the Schwartz kernel of the
resolvent of the Laplacian $\Delta_{\tilde{M}}$, on $\tilde{M}$
the double of  the manifold with boundary $M_0$, and $R_2(x_1,x_2,
\lambda)$ is the Schwartz kernel of the resolvent of
$\Delta_{0,M}$. Using the explicit expressions of the Schwartz
kernels $R_1(x_1,x_2, \lambda)$ and $R_2(x_1,x_2, \lambda)$ and
(\ref{eq:partition resolvent}), it is possible to prove
\begin{lem}\label{lem:S is parametrix}{\upshape \cite[lemma
3.8]{HUS}} For $\lambda \in \C-\R_+$, the operator $S(\lambda)$ is
a parametrix of $R(\lambda)$ in the sense that
$R(\lambda)-S(\lambda)$ is $L^2$-compact.
\end{lem}
The meromorphic extension of the resolvent will have a domain
contained in the following surface, that we call {\it spectral
surface}:
\begin{equation*}
\Sigma_s:=\{\Lambda:=(\Lambda_i)\in \C^\N \text{ : } \forall i, j
\in \N, \text{  } \Lambda_i^2+\mu_i=\Lambda_j^2+\mu_j\}.
\end{equation*}
$\Sigma_s$ is a covering of $\C$ with projection
$\pi_s(\Lambda):=\Lambda_i^2+\mu_i$.

As we said previously, we will consider the resolvent acting on
weighted $L^2$-spaces that we shall define now:
\begin{equation}\label{eq:def weighted L2}
L^2_\delta(M):=\{f:M\stackrel{\text{meas}}{\rightarrow} \C:
\int_0^\infty \int_Y e^{2\delta
u}\normv{f(y,u)}^2dvol_{g_Y}du<\infty \}.
\end{equation}
For all $\delta>0$, we have the inclusions
\begin{equation*}
L^2_\delta(M)\subset L^2(M)\subset L^2_{-\delta}(M).
\end{equation*}
If we define the {\it physical domain} $FD$ by
\begin{equation*}
FD:=\{\Lambda \in \Sigma_s: \Lambda_i \geq 0 \},
\end{equation*}
then we can identify $FD$ with $\C-\R_+$. We denote by
$\Sigma^{\mu}_s$ the connected component of
$\pi_s^{-1}(\C-[\mu,\infty))$ and, for $\epsilon>0$,
$\mu(\epsilon):=\min\{\mu \in \sigma(\Delta_Y):\mu \leq
\epsilon\}$. Now we can define the domains $\Omega_\epsilon$, for
$\epsilon>0$, where we will extend the resolvent of $\Delta_M$,
\begin{equation*}
\Omega_\epsilon:=\left( FD \cup \pi_s^{-1}(\{z \in \C:
\normv{z}\leq \epsilon\}) \right)\cap \Sigma^{\mu(\epsilon)}_s
\end{equation*}
From the explicit formulas of the resolvents $R_1(\lambda)$ and
$R_2(\lambda)$ in equation (\ref{eq:kernel resolv}), and the
definition of the weigthed $L^2$-spaces in (\ref{eq:def weighted
L2}), we deduce the following.
\begin{lem}\label{lem:ext parametrix}{\upshape\cite[lemma 3.20]{HUS}} For all
$\delta \geq \epsilon >0$ the function $\lambda \mapsto
S(\lambda)$ has a meromorphic extension to $\Omega_\epsilon$ as a
continuous operator from $L^2_\delta(M)$ to $L^2_{-\delta}(M)$.
\end{lem}
For $\Lambda \in \Omega_\epsilon$, let us define the operators
\begin{equation*}
G(\Lambda):=S(\Lambda)(\Delta_M-\pi_s(\Lambda))-Id.
\end{equation*}
We will use the following important tool of functional analysis to
meromorphically extend the resolvent.
\begin{thm}\label{thm:analytic Fredholm}{\upshape \cite[appendix]{HUS}
 (Analytic Fredholm theorem)} Let $U
\subset \C$ be an open and connected set and let $T(z)$, for $z
\in U$, be an analytic family of compact operators of a Hilbert
space $\mathscr{H}$. Suppose that, for $z_0 \in U$,
$(Id-T(z_0))^{-1}$ exists, then the family $(Id-T(z))^{-1}$ is
meromorphic in $U$ with values in the bounded linear operators of
$\mathscr{H}$ and poles contained in the set $\{z \in U:1 \in
\sigma (T(z))\}$.
\end{thm}
As for lemma \ref{lem:S is parametrix}, it  is possible to prove
that $G(\Lambda)$, as a bounded operator of $L^2_{-\delta}(M)$, is
compact. This fact, lemma \ref{lem:ext parametrix} and theorem
\ref{thm:analytic Fredholm} imply
 \begin{thm}\label{thm:merom. ext resolv}{\upshape \cite[theorem  3.24]{HUS}} The
resolvent $R(\lambda)$ has a meromorphic extension from $FD$ to
$\Omega_\epsilon$ as a continuous operator from $L^2_\delta(M)$ to
$L^2_{-\delta}(M)$.
\end{thm}
The next theorem provides the connection between the meromorphic
extension of the resolvent and the absence of singular continuous spectrum.
\begin{thm}\label{thm:abssing1}{\upshape \cite{RS3}} Let $H$ be a
 self-adjoint operator with resolvent
 $R(\lambda):=(H-\lambda)^{-1}$.
\begin{itemize}
\item[i)] Let $(a,b)$ be a bounded interval
 and $\varphi \in \mathscr{H}$. Suppose that there exists
 $p>1$ for which:
\begin{equation}\label{eq:estimsingesp}
sup_{0<\epsilon<1}\int_a^b \vert Im(\varphi,R(x+i\epsilon)\varphi)
\vert^pdx<\infty.
\end{equation}
Then $E_{(a,b)}\varphi\in \mathscr{H}_{ac}$.
\item[ii)] Let
$(a,b)$ be a bounded interval. Suppose that there is a dense
subset $D$ in $\mathscr{H}$ so that for $\varphi \in D$ the
inequality (\ref{eq:estimsingesp}) holds for some $p>1$. Then $H$
has purely absolutely continuous spectrum on $(a,b)$.
\end{itemize}
\end{thm}
%------------------------------------------------------------------------------------------
%--------------------------------------------------------
%--------------------------------------------------------
%--------------------------------------------------------
%--------------------------------------------------------
\section{Analytic dilation  on complete manifolds
with corners of codimension 2}\label{sect:Mourreanddilat} In
\cite{MuellerCorner} it has been explained how to meromorphically
 extend the resolvent of a generalized Laplacian on a complete
  manifold with corner of codimension 2,
  using the method outlined in section \ref{section:classic results}
  under the hypothesis that the Laplacian on the corner has
    kernel $\{0\}$. It turns out that to
   weaken this hypothesis and try to  use the methods of section~\ref{section:classic results}
    is not easy and
   new methods should be used to prove absence of singular continuous spectrum.
    In
   this section we survey  the method of analytic dilation
   applied to compatible Laplacians on complete manifolds
   with corners of codimension 2.
    This
    method appeared originally in the context of
    Shr\"odinger operators and was adapted in \cite{CANODILATION}
     to this geometric context.
\\
\\
Following~\cite{MuellerCorner}, we explain the notions of {\it
compact and complete manifolds with corner of codimension 2}. Let
$X_0$ be a compact oriented Riemannian manifold with boundary $M$
and suppose that there exists a hypersurface $Y$ of $M$ that
divides $M$ in two manifolds with boundary $M_1$ and $M_2$, i.e.
$M=M_1 \cup M_2$
  and $Y=M_1 \cap M_2$. Assume also that  a neighborhood of $Y$ in $M$ is diffeomorphic to $Y \times (-\varepsilon,\varepsilon)$. We say that the manifold $X_0$ {\bf has a corner of
codimension~2} if $X_0$ is endowed with a Riemannian metric $g$
that is a
  product metric on small neighborhoods, $M_i \times (-\varepsilon,0]$ of the $M_i$'s and on a small neighborhood, $Y \times (-\varepsilon,0]^2$, of  the corner $Y$. If $X_0$ has a corner of codimension 2, we say that  $X_0$ is a {\bf compact manifold with corner of codimension 2}.
\bigskip
\begin{center}
\psset{unit=0.5cm}
    \begin{pspicture}(4,1)(12,7)
        %\psgrid[subgriddiv=1,griddots=8,gridlabels=8pt](4,1)(12,7)
            \pscurve[](6.5,4)(7,3.7)(7.5,3.5)(8,4)%hueco dona
                \pscurve[](7,3.7)(7.5,3.9)(7.5,3.5)%hueco dona
            \pscurve[](10,6)(7,5.8)(5,5.7)(7,5.5)(9,5.4)%frontera superior
            \pscurve[](10,6)(9.5,3)(9,5.4)%frontera inferiot
            \pscurve[](5,5.7)(6,2.8)(8.5,3)(9.5,3)%frontera "izquierda"
            \rput(7.2,6.5){$M_1$}
            \rput(10.8,4.6){$M_2$}
            \pscircle*(10,6){0.2}
            \pscircle*(9,5.4){0.2}
            \psline(9,5.4)(11.5,6.1)
            \psline(10,6)(11.5,6.1)
            \rput(11.9,6.3){$Y$}
            \rput(8,1.8){Figure 1. Compact manifold with corner of codimension 2.}
        \end{pspicture}
\end{center}
\begin{exam} For $i=1,2$, let $M_i$ be a compact oriented Riemannian manifold with boundary $\partial M_i:=Y_i$. Suppose that on a neighborhood $Y_i \times (-\varepsilon, 0]$ of $Y_i$ the Riemannian metric $g_i$ of $M_i$ is a product metric i.e. $g_i:=g_{Y_i}+du \otimes du$ where $u$ is the coordinate associated to $(-\varepsilon, 0]$ in $Y_i \times (-\varepsilon, 0]$. Then $M_1 \times M_2$ is a compact manifold with corner of codimension 2.\end{exam}
From the compact manifold with corner $X_0$ we construct a
complete manifold $X$. Let
\begin{equation*}
Z_i:=M_i \cup_Y (\R^+ \times Y), \text{ i=1,2},
\end{equation*}
where the bottom $\{0\} \times Y$ of the half-cylinder is
identified with $\partial M_i=Y$. Then $Z_i$ is a
 manifold with cylindrical end (see definition~\ref{defin:manifold cyl end}).
Define the manifolds
\begin{equation*}
W_1:=X_0\cup_{M_2} (\R_+ \times M_2) \text{ and
}W_2:=X_0\cup_{M_1} (\R_+ \times M_1).
\end{equation*}
Observe that $W_i$ is an $n$-dimensional manifold with boundary
$Z_i$ that can be equipped with a Riemannian metric compatible
with the product metric of $\R_+ \times M_2$ and the Riemannian
metric of $X_0$. Set:
\begin{equation*}
X:=W_1 \cup_{Z_1}(\R_+ \times Z_1)=W_2 \cup_{Z_2}(\R_+ \times
Z_2),
\end{equation*}
where we  identify $\{0\}\times Z_i$ with $Z_i$, the boundary of
$W_i$.
\bigskip%%%
\begin{center}
      \psset{unit=0.5cm}
          \begin{pspicture}(0,-1.5)(10,11)
              %\psgrid[subgriddiv=1,griddots=8,gridlabels=8pt](0,0)(11,10)
                  \psline(0,0)(0,10)
                   \qline(0,0)(10,0)
        \psline(2,0)(2,10)
        \psline(0,2)(10,2)
        \rput(1,1){$X_0$}
        \rput(9,9){$[0,\infty)^2 \times Y$}
        \rput(9,2.8){$Z_1$ }
        \rput(2.6,10.5){$Z_2$}
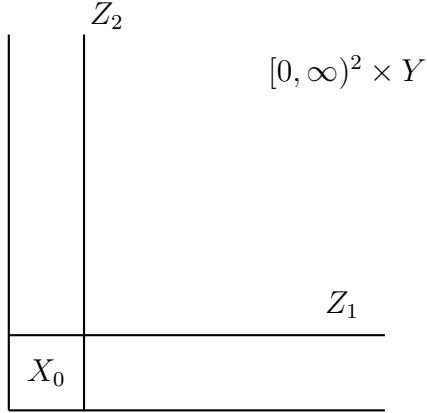
        \rput(7,-1){Figure 2. Sketch of a complete manifold with corner of codimension 2.}
    \end{pspicture}
\end{center}
The above picture is an sketch, in particular the lines that
enclosed the picture should not be thought as boundaries.

Let $T \geq 0$ be given and set
\begin{equation}
Z_{i,T}:=M_i \cup_Y([0,T]\times Y), \text{ for }i=1,2,
\end{equation}
where $\{0\} \times Y$ is identified with $Y$, the boundary of
$M_i$. $Z_{i,T}$ is a family of manifolds with boundary which
exhausts $Z_i$. Next we attach to $X_0$ the manifold $[0,T] \times
M_1$ by identifying $\{0\} \times M_1$ with $M_1$. The resulting
manifold $W_{2,T}$ is compact manifold with corner of codimension
2, whose boundary is the union of $M_1$ and $Z_{2,T}$. The
manifold $X$ has associated a natural exhaustion given by
\begin{equation}\label{eq:def exhaust}
 X_T:=W_{2,T} \cup_{Z_{2,T}}([0,T] \times Z_{2,T}), \text{  }T\geq 0
\end{equation}
where we identify $Z_{2,T}$ with $\{0\} \times Z_{2,T}$.
\bigskip
\begin{center} \psset{unit=0.5cm}
    \begin{pspicture}(0,0)(24,10.5)
        %\psgrid[subgriddiv=1,griddots=8,gridlabels=8pt](0,0)(20,10)
            \pscurve[](6.5,4)(7,3.7)(7.5,3.5)(8,4)%hueco dona
            \pscurve[](7,3.7)(7.5,3.9)(7.5,3.5)%hueco dona
        \pscurve[](5,5.7)(6,2.8)(8.5,3)(9.5,3)
        \pscurve[](10,11)(7,10.8)(5,10.7)(7,10.5)(9,10.4)%cyl. frontera superior
        \psline(9,10.4)(14,10.4)
        \psline(14,10.4)(14,5.3)
        \psline(9.5,3)(14.5,3)
        \psline(5,5.7)(5,10.7)
        \psline(10,11)(15,11)
        \psline(15,11)(15,6)
        \pscurve[](15,6)(14.5,3)(14,5.4)%cyl. frontera inferior
        \rput(9,1){Figure 2. $X_T$, element of the exhaustion of $X$.}
        \psline(9,5.4)(9,10.4)
        \psline(5,5.7)(5,10.7)
        \psline[linestyle=dashed](10,6)(10,11)
        \rput(6.8,8.5){{\tiny $[0,T] \times M_1$}}
 %la linea de abajo la estoy cambiando
        \pscurve[linestyle=dashed](10,6)(8,5.9)(6,5.8)(5.5,5.9)(5,5.7)%frontera superior
\pscurve[](5,5.7)(6,5.4)(7,5.5)(9,5.4)
        \pscurve[linestyle=dashed](10,6)(9.8,3.4)(9.5,3)%frontera inferiot
\pscurve[](9.5,3)(9.2,3.2)(9,5.4)
        \pscurve[](5,5.7)(6,2.8)(8.5,3)(9.5,3)%frontera "izquierda"
        \pscurve[](15,6)(14.5,3)(14,5.4)%cyl. frontera inferior
        \psline(9.5,3)(14.5,3)
        \psline[linestyle=dashed](15,6)(10,6)
        \psline(14,5.4)(9,5.4)
        \rput(11.5,4.7){{\tiny$ [0,T] \times M_2$}}
        \rput(12,8.5){{\tiny $ [0,T]^2\times Y$}}
        \rput(6,4.5){{\tiny $X_0$}}
\end{pspicture}
\end{center}
For each $T \in [0,\infty)$, $X$ has two submanifolds with
cylindrical ends, namely $M_i \times \{T\} \cup (Y \times \{T\})
\times [0,\infty) $, for $i=1,2$.

Let $E$ be a Hermitian vector bundle over a complete manifold with
corner of codimension~2, $X$. Let $\Delta$ be a generalized
Laplacian acting on $C^\infty(X,E)$. The operator $\Delta$ is a
{\bf compatible
  Laplacian} over $X$ if it satisfies the following properties:
\begin{itemize}
\item[i)] There exists a Hermitian vector bundle $E_i$ over $Z_i$
  such that $E \vert_{\R_+\times Z_i}$ is the pullback of $E_i$ under
  the projection $\pi:\R_+\times Z_i \to Z_i$, for $i=1,2$. We suppose
  also that the Hermitian metric of $E$ is the pullback of the
  Hermitian metric of $E_i$. On $\R_+ \times Z_i$, we have
  \begin{equation*}
    \Delta=-\parcial{^2}{u_i^2}+\Delta_{Z_i},
  \end{equation*}
  where $\Delta_{Z_i}$ is a compatible Laplacian acting on
  $C^\infty(Z_i,E_i)$.
\item[ii)] There exists a Hermitian vector bundle $S$ over $Y$
such
  that $E \vert_{\R_+^2 \times Y}$ is the pullback of $S$ under the
  projection $\pi:\R_+^2 \times Y \to Y$. We assume also that the Hermitian product on $E \vert_{\R_+^2 \times Y}$ is the pullback of the Hermitian product on $S$. Finally we suppose that the operator $\Delta$ restricted to $\R_+^2 \times Y$  satisfies
  \begin{equation*}\label{eq:Delta in cylinder}
    \Delta=-\parcial{^2}{u_1^2}-\parcial{^2}{u_2^2}+\Delta_Y,
  \end{equation*}
  where $\Delta_Y$ is a generalized Laplacian acting on
  $C^\infty(Y,S)$.
\end{itemize}
Examples of compatible Laplacians  are the Laplacian acting on
forms and Laplacians associated to compatible Dirac operators (see
 \cite{MuellerCorner}). Since $X$ is a manifold with bounded
geometry and the vector bundle $E$ has bounded Hermitian metric,
the operator $\Delta:C^\infty_c(X,E) \to L^2(X,E)$ is essentially
self-adjoint (see~\cite[corollary 4.2]{Shubin}). Similarly
$\Delta_{Z_i}:C^\infty_c(Z_i,E_i) \to L^2(Z_i,E_i)$ is also
essentially self-adjoint for $i=1,2$.

\begin{defin}
\begin{itemize}
\item Let $H$ and $H^{(i)}$ be the  self-adjoint extensions of
$\Delta:C^\infty_c(X,E) \to L^2(X,E)$ and
$\Delta_{Z_i}:C^\infty_c(Z_i,E_i) \to L^2(Z_i,E_i)$. \item Let
$b_i$ be the self-adjoint extension of
$-\frac{d^2}{du_i^2}:C^\infty_c(\R_+) \to L^2(\R_+)$ obtained by
imposing Dirichlet boundary conditions at $0$. \item Let $H_i$ be
the self-adjoint operator $b_i\otimes Id+Id \otimes H^{(i)}$
acting on $L^2(\R_+)\otimes L^2(Z_i,E_i)$. \item Let $H^{(3)}$ be
the self-adjoint operator associated to the essentially
self-adjoint operator $\Delta_Y:C^\infty(Y,S) \to L^2(Y,S)$ and
let $H_3$ be the self-adjoint operator $H_3:=b_1\otimes Id\otimes
Id+Id\otimes b_2\otimes Id+ Id\otimes Id\otimes H^{(3)}$ acting on
$L^2(\R_+) \otimes L^2(\R_+) \otimes L^2(Y)$. \item The operators
$H_i$ are called {\bf channel operators} for $i=1,2,3$.
    \end{itemize}
\end{defin}

The self-adjoint operators $H_1$ and $H_2$ have a free channel of
dimension 1 (associated to $b_1$ and $b_2$, respectively); the
operator $H_3$ has a free channel of dimension 2 (associated to
$b_1\otimes Id\otimes Id+Id\otimes b_2\otimes Id$).  In some parts
of this text we abuse the notation and denote by $H$, $H_i$, and
$H^{(i)}$
 the
Laplacians acting on distributions and the self-adjoint operators
previously defined.

Along the next section we will use the following notation.  Let
$H$ be a self-adjoint operator acting on a Hilbert space
$\mathscr{H}$. We define the Banach space $\mathscr{H}_2(H)$ as
the domain of $H$ with the norm
$\norm{\varphi}_2:=\norm{(\normv{H}+i)\varphi}$. Similarly, we
define the Banach spaces $\mathscr{H}_1(H)$ as the completion of
$\mathscr{H}_2(H)$ with the norm
$\norm{\varphi}_1:=\norm{(\normv{H}+i)^{1/2}\varphi}$, and
$\mathscr{H}_{-1}(H)$ and $\mathscr{H}_{-2}(H)$ the dual spaces
associated to $\mathscr{H}_1(H)$ and $\mathscr{H}_2(H)$.
%----------------------------------------------------------------
%----------------------------------------------------------------
%-----------------------------------------------------------------
\subsection{Analytic dilation}
\label{section: analytic dilation}The analytic dilation of a
many-body Schr\"odinger operator depends on the analytic dilation
of their subsystem Hamiltonians (see \cite{HS1}). In a similar way
the analytic dilation of $H$ is described in terms of the spectral
theory of the operators $H^{(1)}$, $H^{(2)}$ and $H^{(3)}$,
explained above. For $\theta>0$, we define the operator
$U_{i,\theta}:L^2(Z_i) \to L^2(Z_i)$ that essentially  is the
dilation operator by $\theta+1$ up to a compact set. More
precisely:
\begin{equation*}
U_{i,\theta}f(x)=\begin{cases}
f(x)& \text{ for } x \in M_i.\\
\\
(\theta+1)^{1/2}f((\theta+1)u,y)& \text{ for } x=(u,y) \in [0,\infty) \times Y\\
&\hspace{0.5cm}\text{ and for  } u \text{ big enough},
\end{cases}
\end{equation*}
and $U_{i,\theta}f$ is extended to the whole manifold $Z_i$ in such a way
that it sends $C^\infty_c(Z_i)$ into $C^\infty_c(Z_i)$, and it
becomes a unitary operator on $L^2(Z_i)$.
 We refer to \cite{CANODILATION} for the technical details. Similarly, the operators $U_\theta:L^2(X) \to L^2(X)$ are defined by
\begin{equation*}
U_{\theta}f(x)=\begin{cases}
f(x)& \text{ for } x \in X_0.\\
\\
(\theta+1)^{1/2}U_{i,\theta}f((\theta+1)u_i, z_i)& \text{ for } x=(u_i,z_i) \in [0,\infty) \times  Z_i\\
&\hspace{0.5cm} \text{ and for  } u_i \text{ big enough}.
\end{cases}
\end{equation*}
Again $U_\theta f$ is extended to the whole $X$ in such a way
that, for $f \in C^\infty_c(X)$, $U_\theta f \in C^\infty_c(X)$,
and  $U_\theta$ becomes  a unitary operator in $L^2(X)$.

For $\theta \in [0,\infty)$, define $H_\theta:=U_\theta H
U_\theta^{-1}$, a closed operator with domain
$\mathscr{H}_{2}(H)$. We have  that
\begin{equation*}
\mathscr{H}_{2}(H)=\{f \in L^2(X): \Delta_{dist}f \in L^2(X)\}
\end{equation*}
is the second Sobolev space associated to $(X,g)$. Consider the set:
\begin{equation}\label{def: Gamma}
\begin{split}
\Gamma:=\{\theta:=\theta_0+i\theta_1 \in \C:\theta_0>0, \theta_0
> \normv{\theta_1} \text{ and } Im(\theta)^2<1/2\}.
\end{split}
\end{equation}
We will extend the family $H_\theta$ from $[0,\infty)$ to
$\Gamma$.
\begin{center}
\psset{unit=0.5cm}
    \begin{pspicture}(-3,-2)(12,11)
        %\psgrid[subgriddiv=1,griddots=8,gridlabels=8pt](4,1)(12,7)
                \pspolygon[fillstyle=hlines,linestyle=dotted](12,0)(8,0)(4,4)(8,8)(12,8)
                %\psline(4,4)(8,8)
            %\psline(4,4)(8,0)
            %\psline(8,8)(12,8)
            %\psline(8,0)(12,0)
            \psline{->}(-1,4)(12,4)
            \psline{->}(4,-1)(4,10)
            \rput(3.2,8){$\frac{\sqrt{2}}{2}$}
            \rput(3,0){$-\frac{\sqrt{2}}{2}$}
            \rput(8,3.5){$\frac{\sqrt{2}}{2}$}
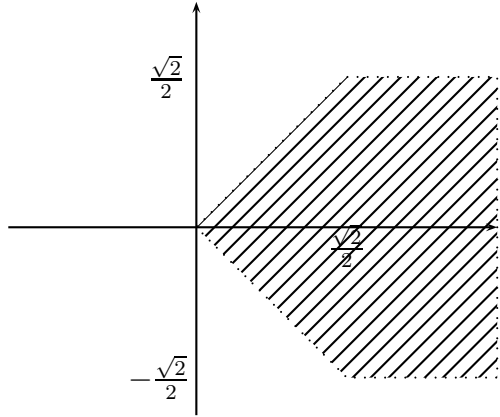
            \rput(4,-1.5){Figure 4. Sketch of the region $\Gamma$.}
        \end{pspicture}
\end{center}
In \cite{CANODILATION}, the next theorem is proved:
\begin{thm}{\upshape \cite{CANODILATION}}\label{thm:Afamily} The family $(H_\theta)_{\theta \in [0,\infty)}$ extends to an holomorphic
family for $\theta \in \Gamma$,  which satisfies:
\\
\\
1) $H_\theta$ is a closed operator with domain
$\mathscr{H}_{2}(H)$ for $\theta \in \Gamma$.
\\
\\
2) For $\varphi \in \mathscr{H}_{2}(H)$ the map $\theta \mapsto
H_\theta \varphi$ is holomorphic in $\Gamma$.
\end{thm}
An holomorphic family of operators  satisfying 1) and 2) is called
{\bf a holomorphic family of type A}. This theorem is proved using
the analogous result that the family $\{H^{(i),\theta}\}_{\theta
\in [0,\infty)}$    extends to a holomorphic family of type A in
$\Gamma$, where $H^{(i),\theta}$ denotes the closed operator
associated to $U_{i,\theta} \Delta_{Z_i} U_{i,\theta}^{-1}$ with
domain
\begin{equation*}
\mathscr{H}_{2}(H^i)=\{f \in L^2(Z_i):\Delta_{dist}(f) \in
L^2(Z_i)\},
\end{equation*}
the second Sobolev space associated to $(Z_i,g_i)$.

The families $H_\theta$ and $H^{(i)}_\theta$ extend to domains
larger than $\Gamma$, but $\Gamma$ is enough for our goals.  In
particular,  $\Gamma$ is chose because for $\theta \in \Gamma$ is
easy to prove that $H^{(i)}_\theta$ is $m$-sectorial (see
\cite[section 2.7]{CANODILATION}), fact that will be important for
the proof of theorem \ref{thm: ess Htheta} where Ichinose lemma is
a main tool.
%-------------------------------------------------------------------------
We define
\begin{equation*}
\theta':=\frac{1}{ (\theta+1)^2}.
\end{equation*}
The parameter $\theta'$ is very important in the description of
the essential spectrum of $H_\theta$ as we can see in  the next
theorem.
\begin{thm}{\upshape \cite{CANODILATION}}\label{thm: ess Htheta}
For $\theta \in \Gamma$,
\begin{equation*}
\begin{split}
\sigma_{ess}(H_\theta)=&\bigcup_{\mu \in \sigma(H^{(3)})} (\mu+\theta'[0,\infty))\\
&\cup  \bigcup_{\lambda_1 \in \sigma_{pp}(H^{(1),\theta})} \left(\lambda_1 +\theta'[0,\infty) \right)\\
&\cup  \bigcup_{\lambda_2 \in \sigma_{pp}(H^{(2,\theta)})}
\left(\lambda_2 +\theta'[0,\infty) \right).
\end{split}
\end{equation*}
\end{thm}
It is possible to associate to $(U_\theta)_{\theta \in
[0,\infty)}$ a set of function $\mathscr{V} \subset
\mathscr{H}_{2}(H)$ that satisfies:
\begin{itemize}
\item[i)]$\mathscr{V}$ is dense in $L^2(X)$. \item[ii)] for
$\varphi \in \mathscr{V}$, $U_\theta \varphi $ is defined for all
$\theta \in \Gamma$. \item[iii)] $U_\theta \mathscr{V}$ is dense
in $L^2(X)$ for all $\theta \in \Gamma$.
\end{itemize}
The elements of a subset of $\mathscr{H}_{2}(H)$  which satisfies
i) and ii) will be called {\bf analytic vectors}.  We denote by
$\Lambda$ the left-hand plane, more explicitly:
\begin{equation}\label{eq: def Lambda}
\Lambda:=\{(x,y) \in \C: x < 0\}.
\end{equation}
We denote by $R(\lambda)$ the resolvent of $H$ and by
$R(\lambda,\theta)$ the resolvent of $H_\theta$.  Using the
general analytic dilation theory of Aguilar-Balslev-Combes (see
\cite{B}) we describe the nature of the spectrum of $H$. This
theory is based on:
\begin{itemize}
\item[i)]  The knowledge of the essential spectrum of $H_\theta$,
provided by theorem  \ref{thm: ess Htheta}. \item[ii)] The
following equation, that is consequence of the unitarity of
$U_\theta$,
\begin{equation}\label{eq:ext resol elem}
\left<R(\lambda)f,g
\right>_{L^2(X)}=\left<R(\lambda,\theta)U_\theta f,U_\theta g
\right>_{L^2(X)},
\end{equation}
for $f,g \in \mathscr{V}$ and $\theta \in [0,\infty)$.
\end{itemize}
Since the right-hand side of (\ref{eq:ext resol elem}) is defined
for $\lambda \in \Lambda$ and $\theta \in \Gamma$, (\ref{eq:ext
resol elem}) provides a meromorphic extension of the functions
$\lambda \mapsto \left<R(\lambda)f,g \right>_{L^2(X)}$ from
$\Lambda$ to $\C-\sigma(H_\theta)$. From this, using
Aguilar-Balslev-Combes, we deduce the following theorem.
\begin{thm}{\upshape \cite{CANODILATION}}
1) For $f,g \in \mathscr{A}$ the function $\lambda \mapsto \langle
R(\lambda)f,g \rangle_{L^2(X)}$ extends from $\Lambda$ to
$\C-\sigma(H_\theta)$.
\\
\\
2) For all $\theta \in \Gamma$, $H_\theta$ has no singular
continuous spectrum.
\\
\\
3) The accumulation points of $\sigma_{pp}(H)$ are contained in
$\{\infty\}\cup \sigma(H^{(3)}) \cup  \cup_{i=1}^2
\sigma_{pp}(H^{(i)})$.
\end{thm}
In the case of manifolds with cylindrical ends we have shown in
section \ref{section:classic results} the absence of singular
continuous spectrum for their Laplacians; in~\cite{DONELLY1},  by giving a
polynomial bound to the growing of the number of
$L^2$-eigenvalues, it is proved that the unique possible
accumulation point of the pure point spectrum of a Laplacian on a
manifold with cylindrical end  is $\infty$.
%////////////////////////////////////////////////////////////////////
%////////////////////////////////////////////////////////////////////%////////////////////////////////////////////////////////////////////
%------------------------------------------------------------------------------------------
%--------------------------------------------------------
%////////////////////////////////////////////////////////////////////
%////////////////////////////////////////////////////////////////////
%////////////////////////////////////////////////////////////////////
%////////////////////////////////////////////////////////////////////
%////////////////////////////////////////////////////////////////////
\section*{Acknowledgment}
The author is very grateful with the organizers of the summer
school of Villa de Leyva 2011 for giving him the opportunity to
participate in this nice school  and to write this paper for its
proceedings. Special thanks are due to Alexander Cardona who
improves a lot the presentation of this paper. He wants to thanks
also Professor Werner M\"uller at Bonn University for his constant
advise that introduced him in the topics of this text and
supported the writing of his PhD-thesis on which section
\ref{sect:Mourreanddilat} and \cite{CANODILATION} are based.
\bibliographystyle{alpha}
\bibliography{literatur}

\end{document}